\documentclass{elsart}

\usepackage{amsmath}
\usepackage{amsfonts}
\usepackage{amssymb}
\usepackage{epsfig}
\usepackage{color}

\makeatletter
\newif\if@restonecol
\makeatother

\newcommand{\ABS}[1]{\left\lvert#1\right\rvert}

\newcommand{\B}{\mathcal{B}}

\newcommand{\ic}{\mathbf{i}} % imaginary unit
\newcommand{\RR}{\mathbb{R}} % real numbers
\newcommand{\CC}{\mathbb{C}} % complex numbers
\newcommand{\TT}{\mathbb{T}} % complex unit circle
\newcommand{\DD}{\mathbb{D}} % open complex unit disc
\renewcommand{\L}{\mathcal{L}} % rational functions
\renewcommand{\P}{\mathcal{P}} % polynomials
\newcommand{\CT}[1]{\overline{#1}}

\newcommand{\Span}{\mathop{\mathrm{span}}\nolimits}

\begin{document}

\begin{frontmatter}
\title{An extension of the associated rational functions on the unit circle}%\thanksref{1}}
%\thanks[1]{The work of the first author is partially supported by the Belgian Network DYSCO (Dynamical Systems, Control, and Optimization),
%funded by the Interuniversity Attraction Poles Programme, initiated by the Belgian State, Science Policy Office.
%The scientific responsibility rests with its authors.}

\author{Karl Deckers\corauthref{1}},%\thanksref{2}},
\corauth[1]{Corresponding author. E-mail: karl.deckers@cs.kuleuven.be.}
%\thanks[2]{This author is a Postdoctoral Fellow of the Research Foundation - Flanders (FWO). E-mail: karl.deckers@cs.kuleuven.be.}
\address{Department of Computer Science, Katholieke Universiteit Leuven, \\ Heverlee, Belgium.}
\author{Mar\'{\i}a Jos\'{e} Cantero, Leandro Moral, and Luis Vel\'{a}zquez}
\address{Department of Applied Mathematics, University of Zaragoza, \\ Zaragoza, Spain.}

\vskip 0.25cm
\address{\sl Dedicated to Professor Franz
Peherstorfer. In memoriam.}
\begin{abstract}
A special class of orthogonal rational functions (ORFs) is presented in this paper. Starting with a sequence of ORFs and the corresponding rational
functions of the second kind, we define a new sequence as a linear combination of the previous ones, the coefficients of this linear combination
being self-reciprocal rational functions. We show that, under very general conditions on the self-reciprocal coefficients, this new sequence
satisfies orthogonality conditions as well as a recurrence relation. Further, we identify the Carath\'{e}odory function of the corresponding
orthogonality measure in terms of such self-reciprocal coefficients.

The new class under study includes the associated rational functions as a particular case. As a consequence of the previous general analysis, we
obtain explicit representations for the associated rational functions of arbitrary order, as well as for the related Carath\'{e}odory function. Such
representations are used to find new properties of the associated rational functions.
\end{abstract}

\end{frontmatter}

\section{Introduction}

Since the fundamental work of Stieltjes and Chebyshev, among others, in the 19$^{\mathrm{th}}$ century, orthogonal polynomials (OPs) have been an
essential tool in the analysis of basic problems in mathematics and engineering. For example, moment problems, numerical quadrature, rational and
polynomial approximation and interpolation, linear algebra, and all the direct or indirect applications of these techniques in engineering are all
indebted to the basic properties of OPs. Mostly orthogonality has been considered on the complex unit circle or on (a subset of) the real line.

Orthogonal rational functions (ORFs) were first introduced by D\v{z}rba\-\v{s}ian in the 1960s. Most of his papers appeared in Russian literature,
but an accessible survey in English can be found in~\cite{Bib:procMMD90,Bib:artKMAB97}. These ORFs are a generalization of OPs in such a way that
they are of increasing degree with a given sequence of poles, and the OPs result if all the poles are at infinity. During the last years, many
classical results of OPs are extended to the case of ORFs.

Several generalizations for ORFs on the complex unit circle and the whole real line have been gathered in book~\cite[Chapt.\ 2--10]{Bib:bBGVHN99}
(e.g.\ the recurrence relation and the Favard theorem, the Christoffel-Darboux relation, properties of the zeros, etc.). Other rational
generalizations can be found in e.g.\ \cite{Bib:artKA07,Bib:artLV08}. Further, we refer to~\cite{Bib:artBC08,Bib:artBCDG08,Bib:artAPEO01} and
to~\cite{Bib:artPA03} for the use of these ORFs in respectively numerical quadrature and system identification, while several results about
matrix-valued ORFs can be found in e.g.\ \cite{Bib:artBBA02,Bib:artBBA03}.

Of course, many of the classical OPs are not defined with respect to a measure on the whole unit circle or the whole real line. Several theoretical
results for ORFs on a subset of the real line can be found in e.g.\ \cite[Chapt.\ 11]{Bib:bBGVHN99} and~\cite{Bib:artKA08,Bib:artKDAB07}. For the
special case in which this subset is a real half-line or an interval, we refer to~\cite{Bib:artABPGEHON03,Bib:artBGHN05}
and~\cite{Bib:artDB07,Bib:artKDJA07,Bib:artKDJAB07,Bib:trJB05,Bib:artJAB03,Bib:artJAB01,Bib:artJDAB03} respectively, while some computational aspects
have been dealt with in e.g.\ \cite{Bib:procKJA06,Bib:artKJA07,Bib:artJA02,Bib:artJAB04,Bib:artJVDAb04,Bib:artJAP07,Bib:artJKAJ08}.

By shifting the recurrence coefficients in the recurrence relation
for OPs and ORFs, the so-called associated polynomials (APs) and
associated rational functions (ARFs) respectively are obtained. ARFs
on a subset of the real line have been studied
in~\cite{Bib:artKDAB08,Bib:procKA08} as a rational generalization of
APs (see e.g.\ \cite{Bib:artWVA91}), while APs on the complex unit
circle, on the other hand, have been studied in~\cite{Bib:artFP96}.
However, so far nothing is known about ARFs on the complex unit
circle, and hence, the main purpose of this paper is to
generalize~\cite{Bib:artFP96} to the case of rational functions,
following the ideas developed by  Professor Franz Peherstorfer.

The outline of this paper is as follows. After giving the necessary theoretical background in Section~\ref{Sec:pre}, in Section~\ref{Sec:fsk} we
recall some basic properties of ORFs on the complex unit circle and their so-called functions of the second kind. Although these properties are
basic, they are partially new in the sense that we prove them in a more general context. Next, in Section~\ref{Sec:nc} we use these ORFs and their
functions of the second kind to define a new class of ORFs on the complex unit circle. The ARFs on the complex unit circle will then turn out to be a
special case of this new class of ORFs, and will be dealt with in Section~\ref{Sec:arf}. We conclude in Section~\ref{Sec:ex} with an example.

\section{Preliminaries}\label{Sec:pre}

The field of complex numbers will be denoted by $\CC$, and for the real line we use the
symbol $\RR$. Let $a\in\CC$, then $\Re\{a\}$ refers to the real part of $a$, while
$\Im\{a\}$ refers to imaginary part. Further, we denote
the imaginary unit by $\ic$. The unit circle and the open unit disc are
denoted respectively by
\[
    \TT = \{z\in\CC:\ABS{z}=1\}\hskip 20pt\mbox{and}\hskip 20pt\DD = \{z\in\CC:\ABS{z}<1\}.
\]
%We will use the blackboard $\NN$ to denote the set of natural numbers.
Whenever the value zero is omitted in the set $X\subseteq\CC$, this will be represented by $X_0$; e.g., $\RR_0=\RR\setminus\{0\}$.

For any complex function $f$, we define the involution operation or substar conjugate by
$f_*(z)=\CT{f(1/\CT{z})}$. With $\P_n$ we denote the space of polynomials of degree less
than or equal to $n$, while $\P$ represents the space of all polynomials.
  Further, the set of complex
functions holomorphic on $X\subseteq\CC$ is denoted by $H(X)$.

Let there be fixed a sequence of complex numbers
$\B=\{\beta_0,\beta_1,\beta_2,\ldots\}\subset\DD$, the rational functions
we then deal with, are of the form
\begin{equation}\label{Eq:rf}
  f_k(z) = \frac{c_k z^k + c_{k-1} z^{k-1} + \dots + c_0}{(1 - \CT{\beta}_1z)
  (1 - \CT{\beta}_2z) \cdots (1 - \CT{\beta}_kz)}, \hskip 25pt k=1,2,\ldots \;.
\end{equation}

The first element $\beta_0$  has no influence in the
rational functions, but it will play a role in the corresponding recurrence. The standard
choice is $\beta_0=0$, but in this paper $\beta_0$ will be free. The reason is that, even if we
choose $\beta_0=0$ for the orthogonal rational functions, the corresponding associated rational functions involve a shift in the poles so
that the related sequence $\{\beta_N, \beta_{N+1},\cdots\}$ starts at some $\beta_N$
which is not necessarily zero.

Note that, whenever $1/\CT{\beta}_k=\infty$ for every $k\geqslant1$, the ``rational functions"
$f_k(z)$ in~\eqref{Eq:rf} are in fact polynomials of degree $k$. Thus the polynomial case is automatically accounted for.

We define the Blaschke\footnote{The factors and products are named after Wilhelm Blaschke,
who introduced these for the first time in~\cite{Bib:artBL15}.} factors for $\B$ as
\begin{equation} \label{ch pre subsec parf eq zkz}
    \zeta_k(z)=\eta_{k}\frac{\varpi_k^*(z)}{\varpi_k(z)},\hskip 20pt
    \eta_{k}=\left\{\begin{array}{ll}\frac{\CT{\beta}_k}{\ABS{\beta_k}},
    & \beta_k\neq0 \\ 1,& \beta_k=0 \end{array}\right.,\hskip 20pt
    k=0,1,2,\ldots,
\end{equation}
where
$$
\varpi_k(z) = 1-\CT{\beta}_kz,\hskip 35pt\varpi_k^*(z) = z\varpi_{k*}(z) = z-\beta_k,
$$
and the corresponding Blaschke products for $\B$ as
\begin{equation} \label{ch pre subsec parf eq bkz}
    B_{-1}(z)=\zeta_0^{-1}(z),\hskip 10pt B_k(z)=B_{k-1}(z)\zeta_k(z),\hskip 15pt
    k=0,1,2,\ldots\; .
\end{equation}
These Blaschke products generate the spaces of rational functions with poles in $1/\CT{\beta}_k$, defined by
$$
\L_{-1}=\{0\}, \; \L_0=\CC,\; \L_n:=\L\{\beta_1,\ldots,\beta_n\}=\Span\{B_0,\ldots,B_n\},
\hskip 15pt n\geqslant1,
$$
and $\L=\cup_{n=0}^{\infty}\L_n$. Let
\begin{equation*}
    \pi_0(z)\equiv1,\quad\pi_k(z)=\prod_{j=1}^k\varpi_j(z),\quad k=1,2,\ldots,
\end{equation*}
then for $k\geqslant1$ we may write equivalently
\begin{equation} \label{ch pre subsec parf eq up}
    B_k(z) = \upsilon_k\frac{\pi_k^*(z)}{\pi_k(z)},\hskip 30pt\upsilon_k =
    \prod_{j=1}^k\eta_{j}\in\TT,
\end{equation}
where $\pi_k^*(z)=z^k\pi_{k*}(z)$, and thus
\[
    \L_n=\{p_n/\pi_n:p_n\in\P_n\},\hskip 25pt n=0,1,2,\ldots\;.
\]
Note that $\L_n$ and $\L$ are rational generalizations of $\P_n$ and $\P$. Indeed, if
$\beta_k=0$ (or equivalently, $1/\CT{\beta}_k=\infty$) for every $k\geqslant0$, the
expression in~\eqref{ch pre subsec parf eq zkz} becomes $\zeta_k(z)=z$ and the expression
in~\eqref{ch pre subsec parf eq bkz} becomes $B_k(z)=z^k$. With the definition of the
substar conjugate we introduce $\L_{n*}=\{f_*:f\in\L_n\}$.

The superstar transformation of a complex function $f_{n}\in\L_n\setminus\L_{n-1}$ is
defined as
\[
    f_n^*(z)=B_n(z)f_{n*}(z).
\]
Note that the factor $B_n(z)$ merely replaces the polynomial with zeros
$\left\{\beta_j\right\}_{j=1}^{n}$ in the denominator of $f_{n*}(z)$ by a polynomial with
zeros $\left\{1/\CT{\beta}_j\right\}_{j=1}^{n}$ so that $\L_n^*:=\{B_n f_*:
f\in\L_n\}=\L_n$. Like in this identity, sometimes we will denote $f^*:= B_n f_*$ when we
only know that $f\in \L_n$, even if $f$ could belong to $\L_k$ for some
$k<n$. At any time, the meaning of the superstar transformation should be clear from the
context.

A complex function $F$ is called a Carath\'{e}odory function (abbreviated C-function) in $\DD$ iff
$$
F\in H(\DD)\qquad\text{and}\qquad\Re\{F(z)\}>0,\;z\in\DD.
$$
Important related functions are the Riesz-Herglotz kernel
$$
D(t,z) = \frac{\zeta_0(t)+\zeta_0(z)}{\zeta_0(t)-\zeta_0(z)} = \frac{\varpi_0^*(t)\varpi_0(z)+\varpi_0^*(z)\varpi_0(t)}{\varpi_0(\beta_0)(t-z)},
$$
and the Poisson kernel
$$
P(t,z) = {1\over 2}\left(D(t,z)+D_*(t,z)\right) =
\frac{\varpi_z(z)\varpi_0(t)\varpi_0^*(t)}{\varpi_0(\beta_0)\varpi_z(t)\varpi_z^*(t)},\hskip
15pt P_n(t) := P(t,\beta_n),
$$
where the substar conjugate is with respect to $t$. Note
that $P_*(t,z)=P(t,z)$ and $P(t,z)=\Re\left\{D(t,z)\right\}\;$ for $\;z\in\TT$.

 To the C-function $F$ we then
associate a hermitian $\left(\mathfrak{L}_F(t^{-k}) =
\overline{\mathfrak{L}_F(t^k)}\right)$ linear functional $\mathfrak{L}_F$ on the set of formal power series
$\sum_{k=-\infty}^\infty c_k t^k$ with complex coefficients, so that
$$
F(z) = \mathfrak{L}_F\{D(t,z)\},
$$
where we understand again that $\mathfrak{L}_F$ acts on
$t$. In the remainder we will assume that $F(\beta_0)=1$, and that the functional
$\mathfrak{L}_F$ is positive definite. Thus, $\mathfrak{L}_F\{1\} = 1$, and for every
$f\in\L$
$$
\mathfrak{L}_F\{f_*\} = \CT{\mathfrak{L}_F\{f\}}\qquad\text{and}\qquad\mathfrak{L}_F\{ff_*\}>0\quad\text{for}\quad f\neq0.
$$
This is equivalent to saying that
$$
\mathfrak{L}_F\{f\} = \int_\TT f(t)\,d\mu(t)
$$
for a positive Borel measure $d\mu$ on the unit circle with $\int_\TT d\mu(t) = 1$.

We say that two rational functions $f,g \in \L$ are
orthogonal with respect to $\mathfrak{L}_F$ ($f\perp_F g$) if
$$
\mathfrak{L}_F\{f g_*\} = 0.
$$
The functions of a sequence $\phi_n\in \L_n\setminus\{0\}$ are called
 orthogonal rational functions (ORFs) if
 $$
 \phi_n \perp_F \L_{n-1}
 $$
 \noindent and they are called orthonormal if at the same time
$$
\mathfrak{L}_F\{\phi_n\phi_{n*}\} = 1.
$$
The orthogonality $\phi_n\perp_F \L_{n-1}$
for a function $\phi_n\in
 \L_n\setminus\{0\}$ ensures that, in fact, $\phi_n\in
 \L_n\setminus \L_{n-1}$.

A sequence of functions $f_n\in\L_n$ is called para-orthogonal when
$f_n\perp_F\L_{n-1}(\beta_n)=\{g\in\L_{n-1}:g(\beta_n)=0\}$. Further, a function
$f_n\in\L_n$ is called $k$-invariant (or, self-reciprocal) iff $f_n^*=kf_n$, $k\in\CC$. Let $\Phi_{n,\tau}$ be given by
\begin{equation}\label{Eq:PO}
\Phi_{n,\tau} = \phi_n+\tau\phi_n^*,\;\;\tau\in\TT.
\end{equation}
Then, it is easily verified that a self-reciprocal rational function is
para-orthogonal exactly when it is proportional to a function with the form~\eqref{Eq:PO}.
 Furthermore, the following theorem has been proved
in~\cite[Thm.\ 5.2.1]{Bib:bBGVHN99}.

\begin{thm}\label{Thm:Q}
The zeros of $\Phi_{n,\tau}$, given by~\eqref{Eq:PO}, are on $\TT$ and they are simple.
\end{thm}

\section{Orthogonal rational functions and functions of the second kind}\label{Sec:fsk}

With the ORFs $\phi_n$ and para-orthogonal rational functions (para-ORFs)
$\Phi_{n,\tau}$ we associate the so-called functions of the second kind:
$$
\psi_n(z) = \mathfrak{L}_F\{D(t,z)[\phi_n(t)-\phi_n(z)]\}+\mathfrak{L}_F\{\phi_n(t)\},\;\; n\geqslant0,
$$
(where we understand that $\mathfrak{L}_F$ acts on
$t$) and
$$
\Psi_{n,\tau} = \psi_n-\tau\psi_n^*,\;\;\tau\in\TT,
$$
respectively. We now have the following two lemmas. The
first one, which is partially stated in~\cite[Lem.\ 4.2.1]{Bib:bBGVHN99}, can be understood as a direct
consequence of the recurrence relation appearing below. The second lemma has been proved in~\cite[Lem.\ 4.2.2]{Bib:bBGVHN99} for $n>0$
(the statement is obvious for $n=0$).\footnote{Although we use a slightly different definition of the
Riesz-Herglotz kernel from the one in~\cite{Bib:bBGVHN99}, the proofs in the reference
remain valid.}

\begin{lem}
The functions $\psi_n$ are in $\L_n\setminus
\L_{n-1}$.
\end{lem}

\begin{lem}
For $n>0$, it holds for every $f\in\L_{(n-1)*}$ and $g\in\zeta_{n*}\L_{(n-1)*}$ that
$$
(\psi_nf)(z) = \mathfrak{L}_F\{D(t,z)[(\phi_nf)(t)-(\phi_nf)(z)]\}+\mathfrak{L}_F\{(\phi_nf)(t)\},
$$
and
$$
-(\psi_n^*g)(z) = \mathfrak{L}_F\{D(t,z)[(\phi_n^*g)(t)-(\phi_n^*g)(z)]\}-\mathfrak{L}_F\{(\phi_n^*g)(t)\}.
$$
The same holds true for $n=0$, when $f,g\in\CC$.
\end{lem}

As in the polynomial case, a recurrence relation and a Favard-type theorem can be derived for ORFs and their functions of the second kind.

\begin{thm}\label{Thm:rec}
The following two statements are equivalent:
\begin{enumerate}
  \item $\phi_n\in \L_n\setminus\{0\}$ and $\phi_n\perp_F\L_{n-1}$, for a certain C-function $F$
  with $F(\beta_0)=1$, and $\psi_n$ is the rational function of the second kind of $\phi_n$.
  \item $\phi_n$ and $\psi_n$ satisfy a recurrence relation of the form
  \begin{multline}\label{Eq:rec}
\left(
  \begin{array}{cc}
    \phi_{n}(z) & \psi_{n}(z) \\
    \phi_{n}^*(z) & -\psi_{n}^*(z) \\
  \end{array}
\right) \\ = u_{n}(z)\left(
                       \begin{array}{cc}
                         1 & \CT{\lambda}_{n} \\
                         \lambda_{n} & 1 \\
                       \end{array}
                     \right)\left(
                                   \begin{array}{cc}
                                     \zeta_{n-1}(z) & 0 \\
                                     0 & 1 \\
                                   \end{array}
                                 \right)\left(
  \begin{array}{cc}
    \phi_{n-1}(z) & \psi_{n-1}(z) \\
    \phi_{n-1}^*(z) & -\psi_{n-1}^*(z) \\
  \end{array}
\right),\hskip 15pt n>0,
\end{multline}
where $\lambda_n\in\DD$, and
\begin{equation}\label{Eq:un}
u_{n}(z) = e_n\left(\begin{array}{cc}
                    \rho_{n} & 0 \\
                    0 & \CT{\rho}_{n}\CT{\eta}_{n-1}\eta_{n}
              \end{array}\right)\frac{\varpi_{n-1}(z)}{\varpi_{n}(z)},\hskip 15pt\ABS{\rho_{n}}=1,\hskip 15pt e_n\in\RR_0,
\end{equation}
and with initial conditions $\phi_{0}=\psi_0\in\CC_0$. \vskip 0.15cm In the special case
of orthonormality, the initial conditions are
$$
\phi_0 = \psi_0 = \varrho,\hskip 25pt \ABS{\varrho}=1,
$$
and the constants $e_n$ are given by
\begin{equation}\label{Eq:en}
e_n^2 = \frac{\varpi_n(\beta_n)}{\varpi_{n-1}(\beta_{n-1})}\cdot\frac{1}{1-\ABS{\lambda_n}^2}.
\end{equation}
\end{enumerate}
\end{thm}
\begin{pf}
$(1)\Rightarrow(2)$ has been proved in~\cite[Thm.\
4.1.1]{Bib:bBGVHN99} and~\cite[Thm.\ 4.2.4]{Bib:bBGVHN99} for
$\phi_n$ and $\psi_n$ respectively, under the assumption
$\beta_0=0$. Further, $(2)\Rightarrow(1)$ has been proved
in~\cite[Thm.\ 8.1.4]{Bib:bBGVHN99}, again under the assumption
$\beta_0=0$. It is easily verified that the proofs in~\cite[Thm.\
4.1.1]{Bib:bBGVHN99} and~\cite[Thm.\ 8.1.4]{Bib:bBGVHN99} remain
valid when $\beta_0\neq0$. Also the proof in~\cite[Thm.\
4.2.4]{Bib:bBGVHN99} where $n>1$ remains valid under the assumption
$\beta_0\neq 0$. So, we only need to prove the recurrence relation
for $\psi_n$ when $n=1$.

First, note that
\begin{equation}\label{Eq:phi1}
\phi_1(t) = \frac{e_1\rho_1}{\varpi_1(t)}[\eta_0\varpi_0^*(t)\phi_0+\CT{\lambda}_1\varpi_0(t)\phi_0^*].
\end{equation}
Thus, from the orthogonality of $\phi_1$, it follows that
\begin{equation}\label{Eq:Lphi1}
\eta_0\phi_0\mathfrak{L}_F\left\{\frac{\varpi_0^*(t)}{\varpi_1(t)}\right\} = -\CT{\lambda}_1\phi_0^*\mathfrak{L}_F\left\{\frac{\varpi_0(t)}{\varpi_1(t)}\right\}.
\end{equation}
From~\eqref{Eq:phi1} together with the definition of $\psi_1$ and $D(t,z)$, we obtain
\begin{multline*}
\psi_1(z) = \mathfrak{L}_F\{D(t,z)[\phi_1(t)-\phi_1(z)]\} \\
\hskip 25pt= \frac{e_1\rho_1\left[\eta_0\varpi_1(\beta_0)\psi_0+\CT{\lambda}_1\CT{\varpi_0^*(\beta_1)}\psi_0^*\right]}{\varpi_1(z)\varpi_0(\beta_0)}\mathfrak{L}_F\left\{\frac{\varpi_0^*(t)\varpi_0(z)+\varpi_0^*(z)\varpi_0(t)}{\varpi_1(t)}\right\} \\
\hskip 37pt= \frac{e_1\rho_1[\eta_0\varpi_0^*(z)\psi_0-\CT{\lambda}_1\varpi_0(z)\psi_0^*]}{\varpi_1(z)}\left[\frac{\varpi_1(\beta_0)}{\varpi_0(\beta_0)}+\frac{\CT{\lambda}_1\CT{\varpi_0^*(\beta_1)}\psi_0^*}{\eta_0\psi_0\varpi_0(\beta_0)}\right]\mathfrak{L}_F\left\{\frac{\varpi_0(t)}{\varpi_1(t)}\right\},
\end{multline*}
where the last equality follows from~\eqref{Eq:Lphi1}. Further, we have that
$$
\varpi_0(t) = \frac{\varpi_0(\beta_0)}{\varpi_1(\beta_0)}\varpi_1(t)+\frac{\CT{\varpi_0^*(\beta_1)}}{\varpi_1(\beta_0)}\varpi_0^*(t).
$$
Consequently,
$$
\mathfrak{L}_F\left\{\frac{\varpi_0(t)}{\varpi_1(t)}\right\} =  \frac{\varpi_0(\beta_0)}{\varpi_1(\beta_0)} - \frac{\CT{\lambda}_1\CT{\varpi_0^*(\beta_1)}\psi_0^*}{\eta_0\varpi_1(\beta_0)\psi_0}\mathfrak{L}_F\left\{\frac{\varpi_0(t)}{\varpi_1(t)}\right\},
$$
so that
$$
\left[\frac{\varpi_1(\beta_0)}{\varpi_0(\beta_0)}+\frac{\CT{\lambda}_1\CT{\varpi_0^*(\beta_1)}\psi_0^*}{\eta_0\psi_0\varpi_0(\beta_0)}\right]\mathfrak{L}_F\left\{\frac{\varpi_0(t)}{\varpi_1(t)}\right\}=1.
$$
\hfill $\Box$
\end{pf}

By means of the recurrence relation in the previous theorem, we obtain the following determinant formula (a similar result has been proved
in~\cite[Cor.\ 4.3.2.(2)]{Bib:bBGVHN99} under the assumption $\beta_0=0$).

\begin{thm}\label{Thm:PB}
Suppose $\phi_n\in\L_n\setminus\{0\}$ and $\phi_n\perp_F\L_{n-1}$, for a certain C-function $F$ with $F(\beta_0)=1$,
 and let $\psi_n\in\L_n\setminus\{0\}$ be the rational function of the second kind of $\phi_n$.
 Then,
\begin{equation}\label{Eq:PB}
    \left(\phi_n^*\psi_n+\phi_n\psi_n^*\right)(z) = d_nP_n(z)B_n(z),\qquad d_n\in\RR_0.
\end{equation}
In the special case of orthonormality, it holds that $d_n=2$.
\end{thm}
\begin{pf}
Since
$$
\left(\phi_0^*\psi_0+\phi_0\psi_0^*\right)(z) \equiv 2\ABS{\phi_0}^2,
$$
the equality in~\eqref{Eq:PB} clearly holds for $n=0$ and
$d_0=2$ in the orthonormal case.

Suppose now that the equality in~\eqref{Eq:PB} holds true for $0\leqslant k<n$\hskip 5pt
 with $d_k=2$ in the orthonormal case. We then continue by
induction for $k=n$. From~\eqref{Eq:rec} it follows that
\begin{multline*}
\left(\phi_n^*\psi_n+\phi_n\psi_n^*\right)(z) \\
\hskip 40pt= e_n^2(1-\ABS{\lambda_n}^2)\frac{\varpi_{n-1}^2(z)}{\varpi_{n}^2(z)}\CT{\eta}_{n-1}\eta_n\zeta_{n-1}(z)\left(\phi_{n-1}^*\psi_{n-1}+\phi_{n-1}\psi_{n-1}^*\right)(z) \\
\hskip 15pt=
e_n^2(1-\ABS{\lambda_n}^2)\frac{\varpi_{n-1}^2(z)}{\varpi_{n}^2(z)}\frac{\CT{\eta}_{n-1}\zeta_{n-1}(z)}
{\CT{\eta}_{n}\zeta_n(z)}\frac{P_{n-1}(z)}{P_{n}(z)}d_{n-1}P_{n}(z)B_{n}(z) \\
= d_nP_{n}(z)B_{n}(z),\hskip 245pt
\end{multline*}
where
\begin{equation}\label{Eq:dn}
d_n = e_n^2
\left[\frac{\varpi_n(\beta_n)}{\varpi_{n-1}(\beta_{n-1})}\frac{1}{1-\ABS{\lambda_n}^2}\right]^{-1}
 d_{n-1},
\end{equation}
\noindent so that $d_n\in\RR_0$ and in the orthonormal case
$d_n=d_{n-1}=2$, due to~\eqref{Eq:en}.
\hfill $\Box$
\end{pf}

Finally, the following interpolation properties hold true for (para-)ORFs and their functions of the second kind.

\begin{thm}\label{Thm:IP}
Suppose that $F$ is a C-function, with $F(\beta_0)=1$, and let $\phi_n$ and $\psi_n$ be in $\L_n\setminus\{0 \}$.
Then the following two statements are equivalent:
\begin{enumerate}
  \item $\phi_n\perp_{F}\L_{n-1}$ and $\psi_n$ is the rational function of the second kind of $\phi_n$.
  \item $\phi_n, \psi_n$ satisfy
  \begin{equation}\label{Eq:IP}
  \left\{\begin{array}{lcl}
  \left(\phi_{n}F+\psi_{n}\right)(z) &=& \zeta_0(z)B_{n-1}(z)g_n(z) \\
  \left(\phi_{n}^*F-\psi_{n}^*\right)(z) &=& \zeta_0(z)B_{n}(z)h_n(z)
  \end{array}\right.,\qquad g_n,h_n\in H(\DD).
  \end{equation}
\end{enumerate}
Besides, the function $g_n$ in (\ref{Eq:IP}) satisfies
$g_n(\beta_n)\neq0$\footnote{From Lemma~\ref{Lem:pos} it will in fact follow that $g_n(z)\neq0$ for every $z\in\DD$.}.
\end{thm}
\begin{pf}
$(1)\Rightarrow(2)$ has been proved in~\cite[Thm.\ 6.1.1]{Bib:bBGVHN99} under the assumption $\beta_0=0$.
 The proof in~\cite[Thm.\ 6.1.1]{Bib:bBGVHN99} remains valid for $\beta_0\neq0$, when replacing $t$ and $z$
 with $\zeta_0(t)$ and $\zeta_0(z)$ respectively. Thus, it remains to prove that the rational functions
 $\phi_n,   \psi_n\in \L_n\setminus\{0\}$ in~\eqref{Eq:IP} are unique up to a common non-zero multiplicative
 factor, as well as the fact that $g_n(\beta_n)\neq0$. We will prove both things simultaneously by induction on $n$.

First, consider the case in which $n=0$. Clearly, $\phi_0,
\psi_0\in\CC_0$ satisfy~\eqref{Eq:IP} iff
$\phi_0=\psi_0$. Furthermore,  $g_0(\beta_0)\neq0$ because,
otherwise, evaluating~\eqref{Eq:IP} at $\beta_0$ would give
$$
\phi_0=-\psi_0,\qquad \qquad\phi_{0}=\psi_0,
$$
hence, $\phi_0=\psi_0=0$, in contradiction with our assumption $\phi_0,\psi_0\in\L_0\setminus\{0\}$.

Next, suppose that for $0\leqslant k<n$ the rational
functions $\phi_{k}$ and $\psi_{k}$ in~\eqref{Eq:IP} are unique up to a non-zero
multiplicative factor, and that $g_{k}(\beta_{k})\neq0$. We then continue by induction to
prove that the same holds true for $k=n$. Let $\tilde{\phi}_n,  \tilde{\varphi}_n \in \L_n\setminus\{0\}$, then
$\tilde{\phi}_n = k_n\phi_n + a_{n-1}$ and $\tilde{\psi}_n=k_n\psi_n + b_n$, with $k_{n}\in\CC$, $a_{n-1}\in\L_{n-1}$ and $b_{n}\in\L_n$.
Assuming
$$
\left\{\begin{array}{lcl}
  \left(\tilde{\phi}_{n}F+\tilde{\psi}_{n}\right)(z) &=& \zeta_0(z)B_{n-1}(z)\tilde{g}_n(z) \\
  \left(\tilde{\phi}_{n}^*F-\tilde{\psi}_{n}^*\right)(z) &=& \zeta_0(z)B_{n}(z)\tilde{h}_n(z)
  \end{array}\right.,\qquad \tilde{g}_n,\tilde{h}_n\in H(\DD),
$$
gives
$$
\left\{\begin{array}{lcl}
  \left(a_{n-1}F+b_{n}\right)(z) &=& \zeta_0(z)B_{n-1}(z)\hat{g}_{n-1}(z) \\
  \left(\zeta_{n}a_{n-1}^*F-b_{n}^*\right)(z) &=& \zeta_0(z)B_{n}(z)h_{n-1}(z)
  \end{array}\right.,\qquad \hat{g}_{n-1},h_{n-1}\in H(\DD).
$$
>From the second equality it follows that $b_n^*$ is of the form $\zeta_nb_{n-1}^*$,
$b_{n-1}\in\L_{n-1}$, and hence, that $b_n=b_{n-1}$. Thus,
$$
\left\{\begin{array}{lcl}
  \left(a_{n-1}F+b_{n-1}\right)(z) &=& \zeta_0(z)B_{n-2}(z)g_{n-1}(z) \\
  \left(a_{n-1}^*F-b_{n-1}^*\right)(z) &=& \zeta_0(z)B_{n-1}(z)h_{n-1}(z)
  \end{array}\right.,\qquad g_{n-1},h_{n-1}\in H(\DD),
$$
with $g_{n-1}=\zeta_{n-1}\hat{g}_{n-1}$. Therefore $a_{n-1},
b_{n-1} \in \L_{n-1}$  are solutions of~\eqref{Eq:IP} for $k=n-1$, but with
$g_{n-1}(\beta_{n-1})=0$. This contradicts the induction hypothesis that
$g_{n-1}(\beta_{n-1}) \neq 0$ unless $a_{n-1}=b_{n-1}=0$ which implies $\tilde{\phi}_n =
k_n \phi_n$, $\tilde{\psi}_n = k_n \psi_n$.

 Finally, let us prove that $g_n(\beta_n) \neq 0$. If
 $g_n(\beta_n)= 0$ then
$$
\left\{\begin{array}{lcl}
  \left(\phi_{n}F+\psi_{n}\right)(z) &=& \zeta_0(z)B_{n}(z)\hat{g}_n(z) \\
  \left(\phi_{n}^*F-\psi_{n}^*\right)(z) &=& \zeta_0(z)B_{n}(z)h_n(z)
  \end{array}\right.,\qquad \hat{g}_n,h_n\in H(\DD).
$$
From~\eqref{Eq:PB} it then follows that
\begin{eqnarray*}
d_nP_n(z)B_n(z) &=& \phi_{n}^*(z)\psi_n(z)+\phi_n(z)\psi_n^*(z) \\
&=& \phi_n^*(z)\left(\phi_{n}(z)F(z)+\psi_{n}(z)\right)-\phi_n(z)\left(\phi_{n}^*(z)F(z)-\psi_{n}^*(z)\right) \\
&=& \zeta_0(z)B_n(z)g(z),\hskip 95pt g\in H(\DD) \\
&=& \zeta_0(z)B_n(z)\frac{\varpi_0(z)p_{n-1}(z)}{\pi_n(z)},\qquad p_{n-1}\in\P_{n-1},
\end{eqnarray*}
where the last equality follows from the fact that $\left(\phi_n^*\psi_n+\phi_n\psi_n^*\right)\in\L_n\cdot\L_n$.
Consequently,
\begin{multline*}
\hskip 25pt\frac{\eta_0\varpi_0^*(z)p_{n-1}(z)}{\varpi_n(z)\pi_{n-1}(z)} \;= \;d_nP_n(z)
\;= \;\frac{\hat{d}_n\varpi_0(z)\varpi_0^*(z)} {\varpi_n(z)\varpi_n^*(z)},\hskip 15pt
\hat{d}_n\in\RR_0 \\ \hskip -50pt\Longrightarrow \hskip 7pt p_{n-1}(z)
= \tilde{d}_n\frac{\varpi_0(z)\pi_{n-1}(z)}{\varpi_n^*(z)}\notin\P_{n-1},\hskip 15pt
\tilde{d}_n\in\CC_0,
\end{multline*}
which contradicts the assumption $p_{n-1}\in\P_{n-1}$.
\hfill $\Box$
\end{pf}

The following theorem directly follows from Theorem~\ref{Thm:IP}, and the definition of $\Phi_{n,\tau}$ and $\Psi_{n,\tau}$.

\begin{thm}\label{Thm:IP2} The para-ORFs $\Phi_{n,\tau}\in\L_n\setminus\{0\}$ and their second kind ones $\Psi_{n,\tau}\in\L_n\setminus\{0\}$ satisfy
\begin{equation}\label{Eq:PP}
  \left\{\begin{array}{lcl}
  \left(\Phi_{n,\tau}F+\Psi_{n,\tau}\right)(z) &=& \zeta_0(z)B_{n-1}(z)g_n(z) \\
  \left(\Phi_{n,\tau}^*F-\Psi_{n,\tau}^*\right)(z) &=& \overline{\tau}\zeta_0(z)B_{n-1}(z)g_n(z)
  \end{array}\right.,\qquad g_n\in H(\DD),
  \end{equation}
  with $g_n(z)\neq0$ for every $z\in\DD$.
\end{thm}
\begin{pf}
The equalities in~\eqref{Eq:PP} have been proved in~\cite[Cor.\ 6.1.2]{Bib:bBGVHN99} under the assumption $\beta_0=0$, but the proof remains valid for $\beta_0\neq0$. So, we only need to prove that $g_n(z)\neq0$ for every $z\in\DD$.

Suppose that there exists $\hat{\beta}_n\in\DD$ such that $g_n(\hat{\beta}_n)=0$. Let us then define $R_n,S_n\in\L\{\beta_1,\ldots,\beta_{n-1},\hat{\beta}_n\}\setminus\{0\}$ as
$$
R_n(z) = \frac{\varpi_n(z)}{\hat{\varpi}_n(z)}\Phi_{n,\tau}(z)\qquad\text{and}\qquad S_n(z) = \frac{\varpi_n(z)}{\hat{\varpi}_n(z)}\Psi_{n,\tau}(z),
$$
where $\hat{\varpi}_n(z)=1-\CT{\hat{\beta}_n}z$. From the first equality in~\eqref{Eq:PP} we obtain that
\begin{eqnarray*}
\left(R_nF+S_n\right)(z) &=& \zeta_0(z)B_{n-1}(z)\frac{\varpi_n(z)}{\hat{\varpi}_n(z)}g_n(z),\qquad g_n\in H(\DD) \\
&=& \zeta_0(z)B_{n-1}(z)\hat{\zeta}_n(z)\varpi_n(z)\frac{g_n(z)}{\hat{\eta}_n\hat{\varpi}_n^*(z)},\qquad\hat{\zeta}_n(z)=\hat{\eta}_n\frac{\hat{\varpi}_n^*(z)}{\hat{\varpi}_n(z)} \\
&=& \zeta_0(z)B_{n-1}(z)\hat{\zeta}_n(z)\tilde{g}_n(z),\qquad \tilde{g}_n\in H(\DD),
\end{eqnarray*}
where the last equality follows from the fact that $g_n(\hat{\beta}_n)=0$. On the other hand,
\begin{eqnarray*}
R_n^*(z) &=& B_{n-1}(z)\hat{\zeta}_n(z)R_{n*}(z) = B_n(z)\frac{\hat{\zeta}_n(z)}{\zeta_n(z)}\frac{\varpi_{n*}(z)}{\hat{\varpi}_{n*}(z)}\Phi_{(n,\tau)*}(z) \\
&=& \frac{\hat{\zeta}_n(z)}{\zeta_n(z)}\frac{\varpi_{n}^*(z)}{\hat{\varpi}_{n}^*(z)}\Phi_{n,\tau}^*(z) = \frac{\hat{\eta}_n}{\eta_n}\frac{\varpi_{n}(z)}{\hat{\varpi}_{n}(z)}\overline{\tau}\Phi_{n,\tau}(z) = \frac{\hat{\eta}_n}{\eta_n}\overline{\tau}R_n(z),
\end{eqnarray*}
and similarly,
$$
-S_n^*(z) = \frac{\hat{\eta}_n}{\eta_n}\overline{\tau}S_n(z).
$$
Consequently,
\begin{equation}\label{Eq:imp}
  \left\{\begin{array}{lcl}
  \left(R_nF+S_n\right)(z) = \zeta_0(z)B_{n-1}(z)\hat{\zeta}_n(z)\tilde{g}_n(z) \\
  \left(R_n^*F-S_n^*\right)(z) = \frac{\hat{\eta}_n}{\eta_n}\overline{\tau}\zeta_0(z)B_{n-1}(z)\hat{\zeta}_n(z)\tilde{g}_n(z)
  \end{array}\right.,\qquad \tilde{g}_n\in H(\DD).
\end{equation}
Now, consider the ORF $\hat{\phi}_n\perp_{F}\L_{n-1}$, with $\hat{\phi}_n\in\L\{\beta_1,\ldots,\beta_{n-1},\hat{\beta}_n\}\setminus\{0\}$, and let $\hat{\psi}_n\in\L\{\beta_1,\ldots,\beta_{n-1},\hat{\beta}_n\}\setminus\{0\}$ denote
the rational function of the second kind of $\hat{\phi}_n$. Theorem~\ref{Thm:IP} states that $\hat{\phi}_n$ and $\hat{\psi}_n$
are (up to a multiplicative factor) the only non-zero rational functions in $\L\{\beta_1,\ldots,\beta_{n-1},\hat{\beta}_n\}$ satisfying
$$
  \left\{\begin{array}{lcl}
  \left(\hat{\phi}_{n}F+\hat{\psi}_{n}\right)(z) &=& \zeta_0(z)B_{n-1}(z)\hat{g}_n(z) \\
  \left(\hat{\phi}_{n}^*F-\hat{\psi}_{n}^*\right)(z) &=& \zeta_0(z)B_{n-1}(z)\hat{\zeta}_n(z)\hat{h}_n(z)
  \end{array}\right.,\qquad \hat{g}_n,\hat{h}_n\in H(\DD).
$$
Moreover, it holds that $\hat{g}_n(\hat{\beta}_n)\neq0$ for this solution. Therefore, there cannot exist rational functions $R_n,S_n\in\L\{\beta_1,\ldots,\beta_{n-1},\hat{\beta}_n\}\setminus\{0\}$ satisfying~\eqref{Eq:imp}.
\hfill $\Box$
\end{pf}

Theorem 6 is the main result of this section. It is the
rational extension of ~\cite[Thm.\ 2.1]{Bib:artFPRS94}. Its importance relies on the fact
that it provides us with a characterization of ORFs and their second kind ones in terms
only of the C-function $F$. Theorem 6 will be the key tool to study the associated ORFs
and their extensions, analogously to a similar analysis of the polynomial case in
\cite{Bib:artFP96}.

\section{A new class of orthogonal rational functions}\label{Sec:nc}

 Analogously as has been done in~\cite{Bib:artFP96}, we will study a new class of ORFs generated by a given sequence of ORFs.
The rational functions of the new class will satisfy a similar recurrence to that one of
the initial ORFs, but starting at some index $r$ and with shifted poles and (rotated)
parameters. The associated rational functions will be a particular case when the starting index is $r=0$ and there
is no rotation of the parameters.

 To introduce the new class, we need to consider spaces of rational functions based on different sequences
 of complex numbers.

Given the sequences of complex numbers $\B_N=\{\beta_0,\beta_1,\ldots,\beta_N\}\subset\DD$, $\hat{\B}_n=\{\hat{\beta}_0,\hat{\beta}_1,\ldots,\hat{\beta}_n\}\subset\DD$ and $\tilde{\B}_r=\{\tilde{\beta}_0,\tilde{\beta}_1,\ldots,\tilde{\beta}_r\}\subset\DD$, we define the spaces of rational functions
\begin{eqnarray*}
\L_N &:=& \L\{\beta_1,\ldots,\beta_N\} = \Span\{B_0(z),B_1(z)\ldots,B_N(z)\},\quad\L_0=\CC, \\
\hat{\L}_n &:=& \L\{\hat{\beta}_1,\ldots,\hat{\beta}_n\} = \Span\{\hat{B}_0(z),\hat{B}_1(z),\ldots,\hat{B}_n(z)\},\quad\hat{\L}_0=\CC \\
\tilde{\L}_r &:=& \L\{\tilde{\beta}_1,\ldots,\tilde{\beta}_r\} = \Span\{\tilde{B}_0(z),\tilde{B}_1(z),\ldots,\tilde{B}_r(z)\},\quad\tilde{\L}_0=\CC,
\end{eqnarray*}
and
\begin{eqnarray*}
\L_{N+n} &:=& \L\{\beta_1,\ldots,\beta_N,\hat{\beta}_1,\ldots,\hat{\beta}_n\} = \L_N\cdot\hat{\L}_n,\;\;N,n\geqslant0, \\
\tilde{\L}_{r+n} &:=& \L\{\tilde{\beta}_1,\ldots,\tilde{\beta}_r,\hat{\beta}_1,\ldots,\hat{\beta}_n\} = \tilde{\L}_r\cdot\hat{\L}_n,\;\;r,n\geqslant0,
\end{eqnarray*}
with the convention that $\L_{n}=\L_{0+n}=\hat{\L}_n=\tilde{\L}_{0+n}=\tilde{\L}_n$,
$$
\L_{N+n-1} = \left\{\begin{array}{ll}
\L_{N+(n-1)} = \L_N\cdot\hat{\L}_{n-1},& \hskip 20pt n>0 \\
\L_{N-1},& \hskip 20pt n=0\;,
\end{array}\right.
$$
and
$$
\tilde{\L}_{r+n-1} = \left\{\begin{array}{ll}
\tilde{\L}_{r+(n-1)} = \tilde{\L}_r\cdot\hat{\L}_{n-1},& \hskip 20pt n>0 \\
\tilde{\L}_{r-1},& \hskip 20pt n=0\;.
\end{array}\right.
$$
Further, we set $\hat{\beta}_0=\beta_N$, and hence, $\hat{\zeta}_0(z) = \zeta_N(z)$ and $\hat{B}_{-1}(z) = 1/\zeta_N(z)$.

The main idea is, starting with ORFs whose poles are defined by
$$
\beta_1,\beta_2,\dots,\beta_N,\hat{\beta}_1,\hat{\beta}_2,\dots,\hat{\beta}_n,
$$
to generate  new rational functions with poles defined by
$$
\tilde{\beta}_1,\tilde{\beta}_2,\dots,\tilde{\beta}_r,\hat{\beta}_1,\hat{\beta}_2,\dots,\hat{\beta}_n.
$$
This is the purpose of the following  theorem.

\begin{thm}\label{Thm:NRFs}
For $N,n,r\geqslant0$, suppose $\phi_{N+n}\in\L_{N+n}\setminus\{0\}$ and
$\phi_{N+n}\perp_F\L_{N+n-1}$, and let $\psi_{N+n}$ denote the rational function of the
second kind of $\phi_{N+n}$. Further, suppose $A$, $B$, $C$ and $D$  are self-reciprocal
rational functions in\break $\L_{N}\cdot\tilde{\L}_{r}$, satisfying the following
conditions:
\begin{equation}\label{Eq:a1}
  \tau_A:=\frac{A^*(z)}{A(z)} = -\frac{B^*(z)}{B(z)}=-\frac{C^*(z)}{C(z)} = \frac{D^*(z)}{D(z)},\qquad \tau_A\in\TT,
\end{equation}
\begin{equation}\label{Eq:a2}
B(\beta_j)\neq0,\qquad j=0,\ldots,N-1,\qquad N>0,
\end{equation}
\vskip-20pt
\begin{equation}\label{Eq:a3}
\left(A-BF\right)(z) =\zeta_0(z)B_{N-1}(z)g(z),\qquad g\in H(\DD),
\end{equation}
and
\begin{equation}\label{Eq:a4}
\left(AD-BC\right)(z) = \zeta_0(z)B_{N-1}(z)f(z),\qquad f\in H(\DD).
\end{equation}
\vskip-7pt
Then the rational functions $G_{r+n}$, $H_{r+n}$, $J_{r+n}$ and $K_{r+n}$, defined by
\begin{multline}\label{Eq:GJ}
\hskip-10pt\left(
            \begin{array}{cc}
              G_{r+n}(z) & J_{r+n}(z) \\
              H_{r+n}(z) & -K_{r+n}(z)
            \end{array}
          \right) =\\
          \hskip -50pt\left(
            \begin{array}{cc}
              \phi_{N+n}(z) & \psi_{N+n}(z) \\
              \phi_{N+n}^*(z) & -\psi_{N+n}^*(z)
            \end{array}
          \right) \left(
            \begin{array}{cc}
              A(z) & C(z) \\
              B(z) & D(z)
            \end{array}
          \right) \left\{c_nP_N(z)B_N(z)\right\}^{-1},\;\;\;c_n\in\RR_0,
\end{multline}
are all in $\tilde{\L}_{r+n}$. Furthermore, $G_{r+n}^*(z)=\tau_AH_{r+n}(z)$ and $J_{r+n}^*(z)=\tau_AK_{r+n}(z)$.
\end{thm}
\begin{pf}
From~\eqref{Eq:GJ},
$$
G_{r+n}(z) ={\phi_{N+n}(z) A(z) + \psi_{N+n}(z) B(z)\over c_n P_N(z) B_N(z)},
$$
and
$$
H_{r+n}(z) ={\phi_{N+n}^*(z) A(z) - \psi_{N+n}^*(z) B(z)\over c_n P_N(z) B_N(z)}.
$$

 Concerning the numerators of $G_{r+n}$ and $H_{r+n}$, \eqref{Eq:IP}
 and~\eqref{Eq:a3} give
\begin{multline}\label{Eq:L1}
\left(\psi_{N+n}^*B-\phi_{N+n}^*A\right)(z) \\
\hskip1pt = \psi_{N+n}^*(z)B(z)-\phi_{N+n}^*(z)\{B(z)F(z)+\zeta_0(z)B_{N-1}(z)g(z)\} \\
\hskip 27pt = -\{\phi_{N+n}^*(z)F(z)-\psi_{N+n}^*(z)\}B(z)-\zeta_0(z)B_{N-1}(z)g(z)\phi_{N+n}^*(z) \\
\hskip 2pt = -\zeta_0(z)B_{N+n}(z)h_{N+n}(z)B(z) - \zeta_0(z)B_{N-1}(z)g(z)\phi_{N+n}^*(z) \\
= \zeta_0(z)B_{N-1}(z)k_1(z),\hskip 40pt k_1\in H(\DD),\hskip 107pt
\end{multline}
and
\begin{multline}\label{Eq:L2}
\left(\phi_{N+n}A+\psi_{N+n}B\right)(z) \\
\hskip 1pt = \phi_{N+n}(z)\{B(z)F(z)+\zeta_0(z)B_{N-1}(z)g(z)\}+\psi_{N+n}(z)B(z) \\
\hskip 18pt = \{\phi_{N+n}(z)F(z)+\psi_{N+n}(z)\}B(z)+\zeta_0(z)B_{N-1}(z)g(z)\phi_{N+n}(z) \\
\hskip 2pt = \zeta_0(z)B_{N+n-1}(z)g_{N+n}(z)B(z) + \zeta_0(z)B_{N-1}(z)g(z)\phi_{N+n}(z) \\
= \zeta_0(z)B_{N-1}(z)k_2(z),\hskip 40pt k_2\in H(\DD).\hskip 107pt
\end{multline}
Since the left hand side of~\eqref{Eq:L1} and~\eqref{Eq:L2} is in $\L_N\cdot\L_N\cdot\tilde{\L}_{r+n}$, it follows that
$$
k_1,k_2\in\left\{\begin{array}{ll}
\L\{\beta_N\}\cdot\L_N\cdot\tilde{\L}_{r+n},&\hskip 30ptN>0 \\
\tilde{\L}_{r+n},&\hskip 30pt N=0\;.
\end{array}\right.
$$
On the other hand, taking the superstar conjugate of~\eqref{Eq:L2}, and using the fact that $A$ and $B$ are self-reciprocal and satisfy~\eqref{Eq:a1}, we obtain that
$$
-\tau_A\left(\psi_{N+n}^*B-\phi_{N+n}^*A\right)(z) = \frac{\zeta_N(z)}{\zeta_0(z)}B_N(z)\hat{B}_n(z)\tilde{B}_r(z)k_{2*}(z),
$$
and hence,
$$
-\tau_A\zeta_0^2(z)k_1(z) = \zeta_N^2(z)\hat{B}_n(z)\tilde{B}_r(z)k_{2*}(z).
$$
Consequently,
$$
k_1,k_2\in\left\{\begin{array}{ll}
\L\{\beta_N,\beta_N\}\cdot\tilde{\L}_{r+n},&\hskip 30ptN>0 \\
\tilde{\L}_{r+n},& \hskip 30ptN=0\;,
\end{array}\right.
$$
and
\begin{multline*}
\hskip-9pt k_1(z) = \frac{\varpi_0^2(z)p_{r+n}(z)}{\varpi_N^2(z)\hat{\pi}_n(z)\tilde{\pi}_r(z)},
\hskip 15pt k_2(z) =
\frac{\varpi_0^2(z)q_{r+n}(z)}{\varpi_N^2(z)\hat{\pi}_n(z)\tilde{\pi}_r(z)}, \hskip 15pt
p_{r+n},q_{r+n}\in\P_{r+n}.
\end{multline*}
Since
$$
c_nP_N(z)B_N(z) = c_n\eta_N\CT{\eta}_0\frac{\varpi_N(\beta_N)}{\varpi_0(\beta_0)}\cdot\frac{\varpi_0^2(z)}{\varpi_N^2(z)}\zeta_0(z)B_{N-1}(z),
$$
it now follows that $G_{r+n}(z)$ and $H_{r+n}(z)$ are in $\tilde{\L}_{r+n}$. Further, we have that
\begin{eqnarray*}
G_{r+n}^*(z) &=& \tilde{B}_r(z)\hat{B}_n(z)G_{(r+n)*}(z) = \frac{\tau_A\left(\phi_{N+n}^*(z)A(z)-\psi_{N+n}^*(z)B(z)\right)}{B_N^2(z)\cdot c_nP_{N*}(z)B_{N*}(z)} \\
&=& \tau_A\frac{\phi_{N+n}^*(z)A(z)-\psi_{N+n}^*(z)B(z)}{c_nP_N(z)B_N(z)} = \tau_AH_{r+n}(z).
\end{eqnarray*}

Finally, proving the statement for $J_{r+n}(z)$ and $J_{r+n}^*(z)=\tau_AK_{r+n}(z)$ can be done in a similar way as before, under the condition that
\begin{equation}\label{Eq:a5}
\left(C-DF\right)(z) = \zeta_0(z)B_{N-1}(z)\hat{g}(z), \qquad \hat{g}\in H(\DD).
\end{equation}
So, it remains to prove that~\eqref{Eq:a5} holds true under the assumptions~\eqref{Eq:a1}--\eqref{Eq:a4}. Clearly, condition~\eqref{Eq:a5} holds true for $N=0$. For $N>0$, it follows from~\eqref{Eq:a3}--\eqref{Eq:a4} that
\begin{multline*}
\{C(z)-D(z)F(z)\}B(z) = \{A(z)-B(z)F(z)\}D(z)-\zeta_0(z)B_{N-1}(z)f(z) \\
\hskip 53pt = \zeta_0(z)B_{N-1}(z)\{g(z)D(z)-f(z)\} \\
= \zeta_0(z)B_{N-1}(z)\tilde{g}(z),\hskip 30pt \tilde{g}\in H(\DD).\hskip 58.5pt
\end{multline*}
Condition~\eqref{Eq:a5} now follows due to assumption~\eqref{Eq:a2}.
\hfill $\Box$
\end{pf}

As a consequence of the previous theorem and Theorem~\ref{Thm:rec}, we have the following corollary.

\begin{cor}\label{Cor:rec}
The rational functions $G_{r+n}$ and $J_{r+n}$, defined as before in Theorem~\ref{Thm:NRFs}, satisfy a recurrence relation of the form
\begin{multline}\label{Eq:rec2}
\left(
  \begin{array}{cc}
    G_{r+n}(z) & J_{r+n}(z) \\
    G_{r+n}^*(z) & -J_{r+n}^*(z) \\
  \end{array}
\right) = \\ v_{r+n}(z) \left(
                       \begin{array}{cc}
                         1 & \CT{\gamma}_{r+n} \\
                         \gamma_{r+n} & 1 \\
                       \end{array}
                     \right) \left(
                                   \begin{array}{cc}
                                     \tilde{\zeta}_{r+n-1}(z) & 0 \\
                                     0 & 1 \\
                                   \end{array}
                                 \right) \left(
  \begin{array}{cc}
    G_{r+n-1}(z) & J_{r+n-1}(z) \\
    G_{r+n-1}^*(z) & -J_{r+n-1}^*(z) \\
  \end{array}
\right),\; n>0,
\end{multline}
where $\gamma_{r+n} = \CT{\tau}_A\lambda_{N+n}$,
\begin{equation}\label{Eq:vn}
v_{r+n}(z)=\frac{c_{n-1}}{c_{n}}u_{N+n}(z),
\end{equation}
and (recall that $\hat{\beta}_0=\beta_N$)
$$
\tilde{\zeta}_{r+n-1}(z) = \zeta_{N+n-1}(z)=\hat{\zeta}_{n-1}(z),
$$
and with initial conditions $G_{r},J_r\in\tilde{\L}_r$.
\end{cor}

Theorem~\ref{Thm:NRFs} provides us with a constructive
method to generate a new class of rational functions starting with a given sequence of
ORFs. As we pointed out before, the new rational functions have the same poles as the
initial ORFs, but with the first $N$ ones substituted by $r$ other poles. Besides,
Corollary~\ref{Cor:rec} states that these new rational functions satisfy a similar
recurrence, but with different initial conditions $G_r$, $J_r$, and shifted and rotated
parameters $\gamma_{r+n} = \bar{\tau}_A \lambda_{N+n}$. Nevertheless, this recurrence
does not guarantee the orthogonality because it depends on the orthogonality of the
initial conditions $G_r$, $J_r$. Our aim is to complete the hypothesis of Theorem~\ref{Thm:NRFs} with
a minimum number of conditions to ensure the orthogonality of the new rational functions.
This is the purpose of the following theorem, which is our main result.

\begin{thm}\label{Thm:NORF}
Let $G_{r+n}(z)$, $J_{r+n}(z)\neq 0 $ be defined as before in
Theorem~\ref{Thm:NRFs}, and suppose  $\tilde{\beta}_r=\beta_N$. Further, assume that the self-reciprocals
$A$, $B$, $C$ and $D$ in $\L_N\cdot\tilde{\L}_r$ satisfy~\eqref{Eq:a1} and~\eqref{Eq:a2}, together with the following conditions:
\begin{equation}\label{Eq:a32}
\left(A-BF\right)(z) =\zeta_0(z)B_{N-1}(z)g(z),\qquad g(z)\in H(\DD),
\end{equation}
with
\begin{equation}\label{Eq:a33}
g(\beta)\neq0\text{ for
}\beta\in\{\tilde{\beta}_0,\tilde{\beta_1},\ldots,\tilde{\beta}_r,\hat{\beta}_1,\ldots,\hat{\beta}_n\},
\end{equation}
\vskip-32pt
\begin{equation}\label{Eq:a42}
\hskip -5pt\left(AD-BC\right)(z) =
\zeta_0(z)B_{N-1}(z)\tilde{\zeta}_0(z)\tilde{B}_{r-1}(z)f(z),\quad f(z)\in
H(\DD)\setminus\{0\},
\end{equation}
and $\tilde{F}$, given by
\begin{equation}\label{Eq:Ft}
\tilde{F}(z) = \frac{-C(z)+D(z)F(z)}{A(z)-B(z)F(z)},
\end{equation}
is a C-function, with $\tilde{F}(\tilde{\beta}_0)=1$. Then $G_{r+n}\perp_{\tilde{F}}\tilde{\L}_{r+n-1}$
(respectively $J_{r+n}\perp_{1/\tilde{F}}\tilde{\L}_{r+n-1}$), and $J_{r+n}$ (respectively $G_{r+n}$)
is the function of the second kind of $G_{r+n}$ with respect to $\tilde{F}$ (respectively, of $J_{r+n}$
with respect to $1/\tilde{F}$).
\end{thm}
\begin{pf}
 Theorem~\ref{Thm:NRFs} implies that $G_{r+n}$, $J_{r+n}\in \tilde{\L}_{r+n}\setminus\{0\}$.
 From~\eqref{Eq:GJ}, ~\eqref{Eq:a42}, ~\eqref{Eq:Ft} and~\eqref{Eq:IP} it follows that
\begin{multline*}
\left\{(A-BF)\left(\tilde{F}G_{r+n}+J_{r+n}\right)\right\}(z) = \frac{\left\{(AD-BC)\left(F\phi_{N+n}+\psi_{N+n}\right)\right\}(z)}{c_nP_N(z)B_N(z)} \\
\hskip -55pt= \tilde{\zeta}_0(z)\tilde{B}_{r-1}(z)\zeta_0(z)B_{N}(z)\hat{B}_{n-1}(z)h(z),
\hskip 25pt h\in H(\DD).
\end{multline*}
This, together with~\eqref{Eq:a32} and the condition on the function $g$, gives
\begin{multline*}
\left(\tilde{F}G_{r+n}+J_{r+n}\right)(z)= \tilde{\zeta}_0(z)\tilde{B}_{r-1}(z)\zeta_{N}(z)\hat{B}_{n-1}(z)\hat{h}(z)\\
\hskip -110pt= \tilde{\zeta}_0(z)\tilde{B}_{r+n-1}(z)\hat{h}(z),\hskip 50pt\hat{h}\in
H(\DD).\hskip 50.5pt
\end{multline*}
Next, assuming that $\tilde{F}$ is a C-function, we also
obtain that
$$
\left(G_{r+n}+\frac{1}{\tilde{F}}J_{r+n}\right)(z)= \tilde{\zeta}_0(z)\tilde{B}_{r+n-1}(z)\tilde{h}(z), \hskip 20pt
\tilde{h}\in H(\DD).
$$

Further, it follows from~\eqref{Eq:GJ}, ~\eqref{Eq:a42},
~\eqref{Eq:Ft} and~\eqref{Eq:IP} that
\begin{multline*}
\CT{\tau}_A\left\{(A-BF)\cdot\left(\tilde{F}G_{r+n}^*-J_{r+n}^*\right)\right\}(z) \\
= \frac{\left\{(AD-BC)\cdot\left(F\phi_{N+n}^*-\psi_{N+n}^*\right)\right\}(z)}{c_nP_N(z)B_N(z)} \\
\hskip -95pt= \tilde{\zeta}_0(z)\tilde{B}_{r-1}(z)\zeta_0(z)B_{N}(z)\hat{B}_{n}(z)h(z),
\hskip 35pt h\in H(\DD).\hskip 10pt
\end{multline*}
This, together with~\eqref{Eq:a32}, the condition on the function $g$, and the assumption that $\tilde{F}$ is a C-function, yields
\begin{multline*}
\left(\tilde{F}G_{r+n}^*-J_{r+n}^*\right)(z)=
\tilde{\zeta}_0(z)\tilde{B}_{r-1}(z)\zeta_{N}(z)\hat{B}_{n}(z)\hat{h}(z)\\
\hskip -120pt=\tilde{\zeta}_0(z)\tilde{B}_{r+n}(z)\hat{h}(z), \hskip 40pt \hat{h}\in
H(\DD),\hskip 70pt
\end{multline*}
and
$$
\left(G_{r+n}^*-\frac{1}{\tilde{F}}J_{r+n}^*\right)(z)=
\tilde{\zeta}_0(z)\tilde{B}_{r+n}(z)\tilde{h}(z), \hskip 20pt \tilde{h}\in H(\DD).
$$
The orthogonality now follows from Theorem~\ref{Thm:IP}.
\hfill  $\Box$
\end{pf}

The orthogonality properties
$G_{r+n}\perp_{\tilde{F}}\tilde{\L}_{r+n-1}$  and
$J_{r+n}\perp_{1/\tilde{F}}\tilde{\L}_{r+n-1}$  imply that the hypothesis of Theorem~\ref{Thm:NORF} ensure that, not only $G_{r+n},J_{r+n}\in \tilde{\L}_{r+n}$, but
$G_{r+n},J_{r+n}\in \tilde{\L}_{r+n}\setminus \tilde{\L}_{r+n-1}$ too.

\begin{rem}
>From Theorem~\ref{Thm:PB} it follows that, under the same conditions as in
Theorem~\ref{Thm:NORF}, it should hold that
$$
\left(G_{r+n}^*J_{r+n}+G_{r+n}J_{r+n}^*\right)(z) = \tilde{d}_n\tilde{P}_{r+n}(z)\tilde{B}_{r+n}(z),
\qquad\tilde{d}_n\in\RR_0.
$$
Indeed, taking the determinant on both sides of~\eqref{Eq:GJ}, we find for $n\geqslant0$ that
\begin{multline}\label{Eq:GJ*}
\left(G_{r+n}^*J_{r+n}+G_{r+n}J_{r+n}^*\right)(z) \\\hskip 20pt=
\CT{\tau}_A\frac{\left\{(AD-BC)\cdot
\left(\phi_{N+n}^*\psi_{N+n}+\phi_{N+n}\psi_{N+n}^*\right)\right\}(z)}{\left[c_nP_N(z)B_N(z)\right]^2} \\
\hskip 49pt= \frac{\tilde{\zeta}_0(z)\tilde{B}_{r-1}(z)\cdot
P_{N+n}(z)\hat{B}_n(z)}{P_N(z)\frac{\varpi_0^2(z)}
{\varpi_N^2(z)}}\hat{f}(z), \hskip 25pt\hat{f}\in H(\DD)\setminus\{0\}\\
 \hskip -125pt= \tilde{P}_{r+n}(z)\tilde{B}_{r+n}(z)\tilde{f}(z),\hskip 170pt
\end{multline}
where
$$
\tilde{f}(z)=\frac{P_{N+n}(z)}{P_N(z)\tilde{P}_{r+n}(z)}\cdot\frac{\tilde{\zeta}_0(z)}{\tilde{\zeta}_r(z)}\cdot\frac{\varpi_N^2(z)}{\varpi_0^2(z)}\cdot\hat{f}(z)\in H(\DD)\setminus\{0\}.
$$
Bearing in mind that the left hand side of ~\eqref{Eq:GJ*} is in
$\tilde{\L}_{r+n}\cdot\tilde{\L}_{r+n}$, it follows that $\tilde{f}\in\left(\L\{\tilde{\beta}_0\}\cdot\tilde{\L}_{r+n-1}\right)\setminus\{0\}$ for $r+n>0$, respectively $\tilde{f}\in\CC_0$ for $r+n=0$. Furthermore, taking
the superstar conjugate of~\eqref{Eq:GJ*}, we obtain that
$$
\left(G_{r+n}^*J_{r+n}+G_{r+n}J_{r+n}^*\right)(z) = \tilde{P}_{r+n}(z)\tilde{B}_{r+n}(z)\tilde{f}_*(z),
$$
and hence,
$$
\tilde{f}(z) = \tilde{f}_*(z) \equiv\tilde{d}_n\in\RR_0.
$$
\end{rem}

\section{Associated rational functions}\label{Sec:arf}

A special class of rational functions, the so-called associated rational functions (ARFs), is obtained when $\tau_A=1$ and $r=0$. ARFs orthogonal on a subset of the real line are investigated in detail in~\cite{Bib:artKDAB08}. Analogously to the case of a subset of the
real line, we define the ARFs on the unit circle as follows.

\begin{defn}\label{Def:ARF}
Suppose that the rational functions $\{\phi_n\}$ and $\{\psi_n\}$, with poles among $\{1/\CT{\beta}_1,\ldots,1/\CT{\beta}_n\}$, satisfy a recurrence relation of the form~\eqref{Eq:rec}. Then, for a given $k\geqslant0$, we call the rational functions $\phi_{n\setminus k}^{(k)}$ and $\psi_{n\setminus k}^{(k)}$ generated by the recurrence formula
\begin{multline*}
\left(
  \begin{array}{cc}
    \phi_{n\setminus k}^{(k)}(z) & \psi_{n\setminus k}^{(k)}(z) \\
    \phi_{n\setminus k}^{(k)*}(z) & -\psi_{n\setminus k}^{(k)*}(z) \\
  \end{array}
\right) \\= u_{n}(z)\left(
                       \begin{array}{cc}
                         1 & \CT{\lambda}_{n} \\
                         \lambda_{n} & 1 \\
                       \end{array}
                     \right)\left(
                                   \begin{array}{cc}
                                     \zeta_{n-1}(z) & 0 \\
                                     0 & 1 \\
                                   \end{array}
                                 \right)\left(
  \begin{array}{cc}
    \phi_{(n-1)\setminus k}^{(k)}(z) & \psi_{(n-1)\setminus k}^{(k)}(z) \\
    \phi_{(n-1)\setminus k}^{(k)*}(z) & -\psi_{(n-1)\setminus k}^{(k)*}(z) \\
  \end{array}
\right),\hskip 25pt n>k,
\end{multline*}
with initial conditions $\phi_{k\setminus k}^{(k)}=\psi_{k\setminus k}^{(k)}\in\CC_0$, the ARFs of $\{\phi_n\}$ and $\{\psi_n\}$ of order $k$.
\end{defn}

Note that the subscript '$\setminus k$' in the definition of the ARFs refers to the fact that the ARFs do not have poles
among $\{1/\CT{\beta}_1,\ldots,1/\CT{\beta}_k\}$. In other words, when shifting the
recurrence coefficients, the poles are shifted too. Defining $\L_{n\setminus k}$ by
$$
\L_{(k-1)\setminus k}=\{0\},\hskip 15pt\L_{k\setminus k}=\CC,\hskip 15pt\L_{n\setminus k}
= \L\{\beta_{k+1},\ldots,\beta_n\},
$$
we have that $ \phi_{n\setminus k}^{(k)},\psi_{n\setminus
k}^{(k)} \in \L_{n\setminus k}\setminus \L_{(n-1)\setminus k}$.

 As an application of Theorems~\ref{Thm:NRFs}
and~\ref{Thm:NORF}, we get an explicit representation of the ARFs and of the function to which they are orthogonal in Theorem~\ref{Thm:ARF}. But first we need the following lemma.

\begin{lem}\label{Lem:pos}
Suppose $\phi_k\in\L_k\setminus\{0\}$ such that $\phi_k\perp_F\L_{k-1}$, and let $\psi_k\in\L_k\setminus\{0\}$ denote the rational function of the
second kind of $\phi_k$. It then holds for every $z\in\DD$ that
$$
 \frac{\ABS{(\phi_kF+\psi_k)(z)}^2-\ABS{(\phi_k^*F-\psi_k^*)(z)}^2}{\ABS{(\Phi_{k,-1}F+\Psi_{k,-1})(z)}^2}>0.
$$
\end{lem}
\begin{pf}
First, note that Theorem~\ref{Thm:IP} implies that
$$
(\Phi_{k,\tau}F+\Psi_{k,\tau})(z)=\zeta_0(z)B_{k-1}(z)[g_k(z)+\tau\zeta_k(z)h_k(z)],\quad g_k+\tau\zeta_kh_k\in H(\DD).
$$
Moreover, from Theorem~\ref{Thm:IP2} it follows that
$$
g_k(z)+\tau\zeta_k(z)h_k(z)\neq 0
$$
for every $\tau\in\TT$ and for every $z\in\DD$. Therefore, we have that either
\begin{equation}\label{Eq:ineq}
G(z):=\frac{\ABS{g_k(z)}^2-\ABS{\zeta_k(z)h_k(z)}^2}{\ABS{g_k(z)-\zeta_k(z)h_k(z)}^2} > 0 \text{ \; for every } z\in\DD,
\end{equation}
or
$$
G(z) < 0 \text{ \; for every } z\in\DD.
$$
However, the second option is not possible because for $z=\beta_k\in\DD$ we get $G(\beta_k)=1>0$ \footnote{From~\eqref{Eq:ineq} it follows that $\ABS{g_k(z)}^2>\ABS{\zeta_k(z)h_k(z)}^2\geqslant0$ for every $z\in\DD$, which proves that in Theorem~\ref{Thm:IP} (2), $g_k(z)\neq0$ for every $z\in\DD$.}. The statement now follows by multiply\-ing the
numerator and denominator in~\eqref{Eq:ineq} with $\ABS{\zeta_0(z)B_{k-1}(z)}^2$. \hfill $\Box$
\end{pf}

\begin{thm}\label{Thm:ARF}
For $n\geqslant k\geqslant0$, suppose $\phi_{n}\in\L_{n}\setminus\{0\}$ and $\phi_{n}\perp_F\L_{n-1}$, and let $\psi_{n}$ denote the rational
function of the second kind of $\phi_{n}$. Then, there exist constants $c_{n,k}\in\RR_0$ such that the ARFs
$\phi_{n\setminus k}^{(k)}$ and $\psi_{n\setminus k}^{(k)}$ are given by
\begin{multline}\label{Eq:ARF}
\left(
            \begin{array}{cc}
              \phi_{n\setminus k}^{(k)}(z) & \psi_{n\setminus k}^{(k)}(z) \\
              \phi_{n\setminus k}^{(k)*}(z) & -\psi_{n\setminus k}^{(k)*}(z)
            \end{array}
          \right) =\\
          \left(
            \begin{array}{cc}
              \phi_{n}(z) & \psi_{n}(z) \\
              \phi_{n}^*(z) & -\psi_{n}^*(z)
            \end{array}
          \right)\left(
            \begin{array}{cc}
              \Psi_{k,-1}(z) & -\Psi_{k,1}(z) \\
              -\Phi_{k,-1}(z) & \Phi_{k,1}(z)
            \end{array}
          \right)\left\{c_{n,k}P_k(z)B_k(z)\right\}^{-1}.\hskip 30pt
\end{multline}
Further, $\phi_{n\setminus k}^{(k)}$ (respectively $\psi_{n\setminus k}^{(k)}$) are orthogonal
with respect to the C-function $F^{(k)}$ (respectively $1/F^{(k)}$), given by
\begin{equation}\label{Eq:Ft2}
F^{(k)}(z) = \frac{\Phi_{k,1}(z)F(z)+\Psi_{k,1}(z)}{\Phi_{k,-1}(z)F(z)+\Psi_{k,-1}(z)},
\end{equation}
with $F^{(k)}(\beta_k)=1$. In the special case in which for every $n\geqslant0$ it holds that
\begin{equation}\label{Eq:cn}
c_{n,k}^2 = d_kd_n,
\end{equation}
where $d_j$ is the constant defined in Theorem~\ref{Thm:PB}, then $\phi_{n\setminus
k}^{(k)}$ and $\psi_{n\setminus k}^{(k)}$ are ortho\-normal.
\end{thm}
\begin{pf}
First, put $\tilde{\beta}_0=\beta_k$. Since $\Phi_{k,\tau}^*=\CT{\tau}\Phi_{k,\tau}$ and
$\Psi_{k,\tau}^*=-\CT{\tau}\Psi_{k,\tau}$, condition~\eqref{Eq:a1} is satisfied by
$A=\Psi_{k,-1}$, $B=-\Phi_{k,-1}$, $C=-\Psi_{k,1}$ and $D=\Phi_{k,1}$  with $\tau_A=1$.
Further, condition~\eqref{Eq:a2} is satisfied too due to Theorem~\ref{Thm:Q}.

Theorems~\ref{Thm:PB} and~\ref{Thm:IP2} imply that
\begin{eqnarray}
\nonumber \left(A-BF\right)(z) &=&
\left(\Psi_{k,-1}+\Phi_{k,-1}F\right)(z) \\
 &=& \zeta_0(z)B_{k-1}(z)g_k(z),\quad g_k\in H(\DD), \label{Eq:T} \\
\nonumber \left(C-DF\right)(z) &=&
-\left(\Psi_{k,1}+\Phi_{k,1}F\right)(z) \\
 &=& \zeta_0(z)B_{k-1}(z)h_k(z),\quad h_k\in H(\DD), \label{Eq:N}
\end{eqnarray}
with $g_k(z)\neq0$ and $h_{k}(z)\neq0$ for every $z\in\DD$, and
\begin{multline*}
\left(\Phi_{k,1}\Psi_{k,-1}-\Phi_{k,-1}\Psi_{k,1}\right)(z) \\= \left\{(\phi_k+\phi_k^*)(\psi_k+\psi_k^*)-(\phi_k-\phi_k^*)(\psi_k-\psi_k^*)\right\}(z) \\
\hskip 5pt= 2\left(\phi_k^*\psi_k+\phi_k\psi_k^*\right)(z) =
 \zeta_0(z)B_{k-1}(z)\frac{\hat{d}_k\varpi_0^2(z)}{\varpi_k^2(z)},\hskip
 25pt\hat{d}_k\in\RR_0.\hskip 3pt
\end{multline*}
Hence, conditions~\eqref{Eq:a32}--\eqref{Eq:a42} are satisfied too. Further, from
\eqref{Eq:ARF} and Theorem~\ref{Thm:PB},
$$
\phi_{k\setminus k}^{(k)} = \psi_{k\setminus k}^{(k)} = \frac{d_k}{c_{k,k}}\neq0.
$$
Consequently,  Corollary~\ref{Cor:rec} and Definition~\ref{Def:ARF} show that the
rational functions $\phi_{n\setminus k}^{(k)}$ and $\psi_{n\setminus k}^{(k)}$, defined
by~\eqref{Eq:ARF}, are the ARFs of order $k$ of $\{\phi_{n}\}$
and $\{\psi_{n}\}$ respectively. Moreover, Theorem~\ref{Thm:rec} ensures that
$\phi_{n\setminus k}^{(k)}$ are orthogonal with respect to a C-function $F^{(k)}$, with
$$
F^{(k)}(\beta_k) = \frac{\psi_{k\setminus k}^{(k)}}{\phi_{k\setminus k}^{(k)}}=1.
$$
Note that, for $\tilde{F}$ given by~\eqref{Eq:Ft}, we also have that
\begin{equation*}%\label{Eq:ftb}
\tilde{F}(\beta_k) = \left\{\frac{(\phi_kF+\psi_k)+(\phi_k^*F-\psi_k^*)}{(\phi_kF+\psi_k)-(\phi_k^*F-\psi_k^*)}\right\}(\beta_k)=1,
\end{equation*}
where the last equality follows from the fact that $(\phi_k^*F-\psi_k^*)(\beta_k)=0$ (see Theorem~\ref{Thm:IP}).
Further, from~\eqref{Eq:T}--\eqref{Eq:N} we get
$$
\tilde{F}(z) = -\frac{h_k(z)}{g_k(z)}\in H(\DD).
$$
Furthermore,
$$
\Re\{\tilde{F}(z)\} = \frac{\ABS{(\phi_kF+\psi_k)(z)}^2-\ABS{(\phi_k^*F-\psi_k^*)(z)}^2}{\ABS{(\Phi_{k,-1}F+\Psi_{k,-1})(z)}^2}>0,\;\;z\in\DD,
$$
due to Lemma~\ref{Lem:pos}. Therefore, $\tilde{F}$ is a C-function,
and hence, the equality for $F^{(k)}$ in~\eqref{Eq:Ft2} follows from Theorem~\ref{Thm:NORF}.

Finally, with $c_{n,k}$ given by~\eqref{Eq:cn}, it holds for $n=k$ that
 $\phi_{k\setminus k}^{(k)} = \psi_{k\setminus k}^{(k)} = 1$, while, for $n>k$, we deduce from~\eqref{Eq:un} and~\eqref{Eq:vn} that
$$
\left[e_{n\setminus k}^{(k)}\right]^2 = \frac{c_{n-1,k}^2}{c_{n,k}^2}e_{n}^2 = \frac{d_{n-1}}{d_{n}}e_{n}^2 = \frac{\varpi_n(\beta_n)}{\varpi_{n-1}(\beta_{n-1})}\cdot\frac{1}{1-\ABS{\lambda_n}^2},
$$
where we have applied~\eqref{Eq:dn} in the last equality. Then, the orthonormality is a consequence of~\eqref{Eq:en} in Theorem~\ref{Thm:rec} $(2)$.
\hfill $\Box$
\end{pf}

Based on the previous theorem, the following relations between ARFs of different order can be proved.

\begin{cor}
For $0\leqslant j\leqslant k\leqslant n$, let $K_{n,k}^{(j)}$ be defined by
$$
K_{n,k}^{(j)} = \frac{d_{k\setminus j}^{(j)}}{c_{(n\setminus j),(k\setminus j)}^{(j)}}.
$$
Then, the following relations hold:
\begin{eqnarray}
(a) &&
2K_{n,k}^{(j)}\phi_{n\setminus j}^{(j)}(z) = \left[\left(\phi_{k\setminus j}^{(j)}+\phi_{k\setminus j}^{(j)*}\right)\phi_{n\setminus k}^{(k)}+\left(\phi_{k\setminus j}^{(j)}-\phi_{k\setminus j}^{(j)*}\right)\psi_{n\setminus k}^{(k)}\right](z).\label{Eq:rel1} \\
(b) &&
2K_{n,k}^{(j)}\phi_{n\setminus j}^{(j)}(z) = \left[\left(\phi_{n\setminus k}^{(k)}+\psi_{n\setminus k}^{(k)}\right)\phi_{k\setminus j}^{(j)}+ \left(\phi_{n\setminus k}^{(k)}-\psi_{n\setminus k}^{(k)}\right)\phi_{k\setminus j}^{(j)*}\right](z).\label{Eq:rel2} \\
\nonumber (c) &&
2\frac{d_{n\setminus j}^{(j)}}{d_{k\setminus j}^{(j)}}K_{n,k}^{(j)}P_{n\setminus k}(z)B_{n\setminus k}(z)\phi_{k\setminus j}^{(j)}(z) \\
&&\qquad\qquad= \left[\left(\psi_{n\setminus k}^{(k)*}+\phi_{n\setminus k}^{(k)*}\right)\phi_{n\setminus j}^{(j)}+\left(\psi_{n\setminus k}^{(k)}-\phi_{n\setminus k}^{(k)}\right)\phi_{n\setminus j}^{(j)*}\right](z), \label{Eq:rel3} \\
\nonumber&&\text{where}\quad P_{n\setminus k}(z) = \frac{P_n(z)}{P_k(z)}\quad\text{and}\quad B_{n\setminus k}(z) = \frac{B_n(z)}{B_k(z)}.
\end{eqnarray}
In the special case in which all the involved ARFs  are
orthonormal, it holds that $K_{n,k}^{(j)}=1$. Further, the corresponding relations for
$\psi_{n\setminus j}^{(j)}$ are obtained by replacing $\phi$ with $\psi$
in~\eqref{Eq:rel1}--\eqref{Eq:rel3}, and vice versa.
\end{cor}
\begin{pf}
It suffices to prove the relations for $j=0$. Relation~\eqref{Eq:rel1} follows immediately from~\eqref{Eq:ARF},
 with the help of the identity
\begin{multline*}
\left[(\Phi_{k,1}\Psi_{k,-1}-\Phi_{k,-1}\Psi_{k,1})\right](z) \\= \left[(\phi_k+\phi_k^*)(\psi_k+\psi_k^*)-(\phi_k-\phi_k^*)(\psi_k-\psi_k^*)\right](z) = 2d_kP_k(z)B_k(z).
\end{multline*}
Next, note that~\eqref{Eq:rel2} is just a reformulation of~\eqref{Eq:rel1}.

Finally, from~\eqref{Eq:ARF} and~\eqref{Eq:PB}, we get
\begin{eqnarray*}
  \left(\phi_{n\setminus k}^{(k)}\phi_n^*-\phi_{n\setminus k}^{(k)*}\phi_n\right)(z) &=& \frac{d_n}{c_{n,k}}P_{n\setminus k}(z)B_{n\setminus k}(z)\left[\phi_k^*(z)-\phi_k(z)\right] \\
  \left(\psi_{n\setminus k}^{(k)}\phi_n^*+\psi_{n\setminus k}^{(k)*}\phi_n\right)(z) &=& \frac{d_n}{c_{n,k}}P_{n\setminus k}(z)B_{n\setminus k}(z)\left[\phi_k(z)+\phi_k^*(z)\right].
\end{eqnarray*}
Relation~\eqref{Eq:rel3} now follows immediately by substraction.
\hfill $\Box$
\end{pf}

\section{Example}\label{Sec:ex}
  In this section we will illustrate the preceding results with an example.
  We will consider the orthonormal rational functions with respect to the Lebesgue measure
  $$d\mu(z) = {dz\over 2\pi\ic z} = {d\theta \over 2\pi},\hskip 25pt z=e^{\ic\theta},
  $$
  and poles defined by $\beta_0, \beta_1,\beta_2,\beta_3,\dots$
   with $\beta_0=0$.
  It is very well known that in this case the parameters $\lambda_n$ of  recurrence~\eqref{Eq:rec} vanish for all $n>0$, so that
 $$
 \phi_n(z) = \sqrt{\omega_n(\beta_n)} \; {z\over \omega_n^*(z)}\; B_n(z) =
 \psi_n(z), \hskip 20pt \phi_n^*(z) = \sqrt{\omega_n(\beta_n)}\; {1\over \omega_n(z)} =
 \psi_n^*(z),
 $$
 where we have used the notation of the previous sections. The corresponding Carath\'{e}odory
 function is
 $$
 F(z) = \int_{\TT} {\zeta_0(t) + \zeta_0(z)\over \zeta_0(t) - \zeta_0(z)} \, d\mu(t) =\int_{\TT}{t+z\over t-z} \, d\mu(t) = 1.
 $$
  The ARFs of order 1 are obtained when $r=0$, $\tau_A=1$, $N=1$, and $\tilde{\beta}_0=\beta_1$, thus the related
  poles are defined by  $\beta_1,\beta_2,\beta_3,\dots\;$.
  By Theorem~\ref{Thm:ARF},
  $$
  A(z)=\psi_1^*(z)+\psi_1(z) = \phi_1^*(z)+\phi_1(z) =D(z),$$
  $$ B(z)=\phi_1^*(z)- \phi_1(z) = \psi_1^*(z)-
  \psi_1(z)
  =C(z),
  $$
  so that the orthonormal ARFs of order 1 and the functions of the second kind are given by
  $$\hskip -45pt\phi_{n\setminus 1}^{(1)}(z) =  \psi_{n\setminus 1}^{(1)}(z)={1\over 2 P_1(z)
B_1(z)} \left[A(z)\phi_n(z)
     +B(z)\psi_n(z)\right]
     $$
     $$
     \hskip 52pt={1\over 2 P_1(z) B_1(z)} \left[\left(\phi_1^*(z)+\phi_1(z)\right)
     \phi_n(z) + \left(\phi_1^*(z)-\phi_1(z)\right)\phi_n(z)\right]
     $$
     $$\hskip-131pt=\sqrt{{\omega_n(\beta_n)\over
      \omega_1(\beta_1)}}{\omega_1^*(z)\over \omega_n^*(z)}B_{n\setminus 1}(z).
      $$
 The corresponding  Carath\'{e}odory function is again
$$ F^{(1)}(z) = {-C(z)+D(z) F(z)\over A(z)-B(z) F(z)} =1,$$
but the orthogonality measure $d\mu^{(1)}$ is not the Lebesgue measure because it must
satisfy now
$$ F^{(1)}(z)=\int_{\TT} {\zeta_1(t) + \zeta_1(z)\over \zeta_1(t) - \zeta_1(z)}
\,d\mu^{(1)}(t) = 1.$$ Taking into account that $${\zeta_1(t) + \zeta_1(z)\over \zeta_1(t)
- \zeta_1(z)} = {\omega_1^*(t) \omega_1(z) + \omega_1^*(z) \omega_1(t)\over
\omega_1(\beta_1) (t-z)}$$ it is easy to see that $$ d\mu^{(1)}(z) =
{\omega_1(\beta_1)\over \omega_1(z)\omega_1^*(z)} {dz\over 2\pi \ic} =
{\omega_1(\beta_1)\over |e^{\ic\theta}-\beta_1|^2}{d\theta\over 2\pi}, \hskip 25pt
z=e^{\ic\theta},$$ which is a rational modification of the Lebesgue measure.

\section*{Acknowledgements}

The first author is a Postdoctoral Fellow of the Research Foundation - Flanders (FWO). The work of this author is partially supported by the Belgian Network DYSCO (Dynamical Systems, Control, and Optimization), funded by the Interuniversity Attraction Poles Programme, initiated by the Belgian State, Science Policy Office.
The scientific responsibility rests with its authors.

Most of this work was done during a recent stay of the first author at the University of Zaragoza, Spain. This author is very grateful to Professor Mar\'{\i}a Jos\'{e} Cantero, Professor Leandro Moral, and Professor Luis V\'{e}lazquez, for their hospitality at the University of Zaragoza, and for the valuable discussion during the development of this work.

The work of the second, third and fourth authors is partially supported by a research
grant from the Ministry of  Science and Innovation of Spain, project code
MTM2008-06689-C02, and by Project E-64 of Diputac\'{\i}on General de Arag\'{o}n (Spain).

\end{document}